# Some Solutions of Fractional Order Partial Differential Equations Using Adomian Decomposition Method

Iqra Javed[a], Ashfaq Ahmad[b,c], Muzammil Hussain[d], S. Iqbal[a].

[a] Department of Informatics and System, University of Management and Technology, Lahore, Pakistan

[b] Department of Electronics, The University of Lahore, Lahore, Pakistan
[c] School of Engineering Science, University of Science and Technology of China, China

[d] Department of Computer Science, University of Management and Technology, Lahore, Pakistan

## Abstract

The Adomian decomposition method is a semi-analytical method for solving ordinary and partial nonlinear differential equations. The aim of this paper is to apply Adomian decomposition method to obtain approximate solutions of nonlinear fractional order partial differential equations with fractional derivatives. The fractional derivatives are taken in the sense of Caputo. The solutions of fractional PDEs are calculated in the form of convergent series. Approximate solutions obtained through the decomposition method have been numerically evaluated, and presented in the form of graphs and tables, and then these solutions are compared with the exact solutions and the results rendering the explicitness, effectiveness and good accuracy of the applied method. Finally, it is observed that the applied method (i.e. Adomian decomposition method) is prevailing and convergent method for the solutions of nonlinear fractional-order partial deferential problems.

## 1. Introduction

In the past decade, mathematicians have devoted effort to the study of explicit and numerical solutions to nonlinear fraction differential equations [1, 2]. A number of methods have been presented [9-12], however, Adomian decomposition method [10, 11] gained tremendous popularity and provides rather effective procedure for explicit and numerical solutions of a broad class of differential systems representing real physical problems. An extensive amount of research has been done on fractional calculus, such as [3-8]. The recent research and work related to fractional calculus, particularly on fractional differential equations, can be found in the book of Podlubny [9]. Over the past few years, a number of fractional calculus applications are being used and in the field of science, engineering and economics [4]. Research on non-linear partial differential equations and linearization techniques has gained quite a momentum due the rapidly proliferating use and recent developments of fractional calculus in these fields [10]. To solve equations of different categories like linear or nonlinear, ordinary differential or partial differential equations, integer or fractional etc, a number of methods have been used, such as Adomian's decomposition method [11-15], homotopy perturbation method [16,17], He's variational iteration method [18-20], homotopy analysis method [21], Galerkin method [22], collocation method [23] and other methods [24-26].

Among a multitude of available techniques to tackle nonlinear equations, Adomian Decomposition Method (ADM) has gained tremendous popularity. Adomian's decomposition method introduced by Mathematician George Adomian in 1984 [27-33], is a powerful tool to solve large amount of real world problems [34-38]. The decomposition method is an effective technique to get analytic solutions of the extensive class of dynamical systems without closure approximation, perturbation theory, assumptions of linearization or weak nonlinearity or restrictive assumptions on stochasticity. Although perturbation methods only give the approximate solutions of non-linear fractional differential equations. Generally, there is no such method that gives exact solution. Adomian's decomposition method also provides an approximate solution of nonlinear equations by keeping the problem into its original form.

The aim of this research is to get the numerical solution of the nonlinear fractional partial differential equation using ADM with time-space-fractional derivative of the form

$$\frac{\partial^\alpha u(x,y)}{\partial y^\alpha} + u(x,y)\frac{\partial^\beta u(x,y)}{\partial x^\beta} = g(x) \qquad 0 < \alpha < 1,\ 0 < \beta < 1 \qquad (1.1)$$

Where $u(x,y)$ is assumed to be a function of space, i.e. disappearing for $x < 0$ and $y < 0$ and $\alpha$ and $\beta$ are parameters to describe the order of the fractional derivatives of y and x respectively. The Caputo sense is used while considering fractional derivatives [39]. If we give values of $\alpha = 1$ or $\beta = 1$, the above fractional differential equations will reduce to the standard non-linear partial differential equation.

Finally, this article is organized as: Section 2 is devoted for some basic definitions on the subject of fractional calculus. The details of the ADM formulation for fractional partial differential equation are discussed in Section 3. Section 4 presents some numerical examples with their solutions. In the end, Section 5 concludes this research article.

## 2. Preliminary

In this section, some basic definitions of fractional calculus are recalled that are used in the analysis given in section 3 and 4. It is well-known that there are different definitions of fractional Integral and fractional derivatives, such as, Grünwald-Letnikov, Riesz, Riemann-Liouville, Caputo, Hadamard and Erdélyi-Kober etc [37-40]. However, this paper gives attention to two of them which are Reimann-louville integral operator and Caputo fractional derivative.

**Definition 1:** The Riemann–Liouville fractional integral operator of order $\alpha \geq 0$ of a function $f \in C_\mu, \mu \geq -1$ is defined as

$$J^\alpha f(x) = \frac{1}{\Gamma(\alpha)} \int_a^x (x-\mu)^{\alpha-1} f(\mu) d\mu, \qquad \alpha > 0, x > 0 \qquad (2.1)$$

$$J^0 f(x) = f(x) \qquad (2.2)$$

When we formulate the model of real world problems with fractional calculus, The Riemann–Liouville have certain disadvantages. Caputo proposed in his work on the theory of viscoelasticity [36] a modified fractional differential operator $D_*^\alpha$.

**Definition 2:** The fractional derivative of $f(x)$ in Caputo sense is defined as

$$D^\alpha f(x) = J^{m-\alpha} D^m f(x) = \frac{1}{\Gamma(m-\alpha)} \int_0^x (x-\eta)^{m-\alpha-1} f^{(m)}(\eta) d\eta$$

For

$$m - 1 < \alpha \leq m, \ m \in N, \ x > 0, f \in C_{-1}^m \tag{2.3}$$

**Definition 3:** The Caputo space-fractional derivative operator of order $\alpha > 0$ is defined as

$$D_y^\alpha u(x,y) = \frac{\partial^\alpha u(x,y)}{\partial y^\alpha} = \frac{1}{\Gamma(m-\alpha)} \int_0^y (y-\tau)^{m-\alpha-1} \frac{\partial^m u(x,\tau)}{\partial \tau^m} d\tau, \quad \text{if } m-1 < \alpha \leq m, \ m \in N \tag{2.4}$$

And the Caputo fractional derivative operator of order $\beta > 0$ is defined as

$$D_x^\beta u(x,y) = \frac{\partial^\beta u(x,y)}{\partial x^\beta} = \frac{1}{\Gamma(m-\alpha)} \int_0^x (x-\eta)^{m-\alpha-1} \frac{\partial^m u(\eta,y)}{\partial \eta^m} d\eta, \quad \text{if } m-1 < \alpha \leq m, \ m \in N \tag{2.5}$$

**Definition 4:** If $m - 1 < \alpha \leq m, \ m \in N, \ f \in C_{-1}^m, \mu \geq -1$, then

$$D^\alpha J^\alpha f(x) = f(x)$$

and

$$D^\alpha J^\alpha f(x) = f(x) - \sum_{k=0}^{m-1} f^{(k)}(0^+) \frac{(x-a)^k}{k!}, \ x > 0. \tag{2.6}$$

One can found the properties of the operator $J^\alpha$ in [43], we mention the following:

For $f \in C_\mu^m, \alpha, \beta > 0, \mu \geq -1$ and $\gamma \geq -1$

- $J^\alpha f(x)$ exists for every $x \in [a,b]$
- $J^\alpha J^\beta f(x) = J^{\alpha+\beta} f(x)$
- $J^\alpha J^\beta f(x) = J^\beta J^\alpha f(x)$
- $J^\alpha (x-a)^\gamma = \frac{\Gamma(\gamma+1)}{\Gamma(\alpha+\gamma+1)} (x-a)^{\alpha+\gamma}.$

## 3. ADM Formulation for Fractional PDE's

The Adomian's decomposition method (ADM) can be briefly studied in [7-11]. The formulation using ADM for fractional order PDE's is given in the following steps:

(a) Write the governing fractional order partial differential equation as:

$$\frac{\partial^\alpha u(x,y)}{\partial y^\alpha} + u(x,y) \frac{\partial^\beta u(x,y)}{\partial x^\beta} = g(x) \quad 0 < \alpha < 1, \ 0 < \beta < 1, \tag{3.1}$$

subject to the initial conditions

$$u(x, 0) = f(x), \qquad (3.2)$$

where $\frac{\partial^\alpha}{\partial y^\alpha}$ and $\frac{\partial^\beta}{\partial x^\beta}$ denotes the Caputo fractional derivative operator (2.5).

(b) Construct the decomposition method for fractional order partial differential equation (3.1). It can be written in operators form as:

$$D_y^\alpha\big(u(x,y)\big) + Nu = g(x), \qquad (3.3)$$

where $Nu = u(x,y)D_x^\beta\big(u(x,y)\big)$ is a nonlinear term and $D_y^\alpha = \frac{\partial^\alpha}{\partial y^\alpha}$ and $D_x^\beta = \frac{\partial^\beta}{\partial x^\beta}$ are the Caputo fractional derivative operators respectively. It is clear that $J_y^\alpha$ and $J_x^\beta$ are the Riemann–Liouville fractional integral operator and inverse of $D_y^\alpha$ and $D_x^\beta$ respectively.

(c) Rearrange the equation (3.3) and apply the inverse operator $J_y^\alpha$ to get

$$J_y^\alpha D_y^\alpha\big(u(x,y)\big) = J_y^\alpha(g(x)) - J_y^\alpha(Nu) \qquad (3.4)$$

From which it follows that

$$u(x,y) = \sum_{k=0}^{m-1} \frac{\partial^k}{\partial y^k} u(x, 0^+) \frac{y^k}{k!} + J_y^\alpha(g(x)) - J_y^\alpha(Nu) \qquad (3.5)$$

The decomposition method consist of decomposing the unknown function $u(x, y)$ into a sum of components defined by the decomposition series

$$u(x,y) = \sum_{n=0}^{\infty} u_n(x,y) \qquad (3.6)$$

(d) Expand the nonlinear term $Nu = uu_x$ using Adomian's polynomials $(A_n)$ as:

$Nu = \sum_{n=0}^{\infty} A_n$, These Adomian polynomials can be calculated for all types of nonlinearity by using algorithm proposed by Adomian [25]. In this case, we use the general formula for Adomian polynomials $A_n$ as:

$$A_n = \frac{1}{n!}\left[\frac{\partial^n}{\partial \lambda^n}\{N(\sum_{i=0}^{\infty} \lambda^i u_i)\}\right]_{\lambda=0} \qquad (3.7)$$

All terms of Adomian polynomials can be derived from this formula. The first few terms extracted from this formula are given below:

$$A_0 = u_0(x,y)\frac{\partial^\beta}{\partial x^\beta} u_0(x,y) \qquad (3.8)$$

$$A_1 = u_0(x,y)\frac{\partial^\beta}{\partial x^\beta}u_1(x,y) + u_1(x,y)\frac{\partial^\beta}{\partial x^\beta}u_0(x,y) \qquad (3.9)$$

$$A_2 = u_0(x,y)\frac{\partial^\beta}{\partial x^\beta}u_2(x,y) + u_1(x,y)\frac{\partial^\beta}{\partial x^\beta}u_1(x,y) + u_2(x,y)\frac{\partial^\beta}{\partial x^\beta}u_0(x,y) \qquad (3.10)$$

(e) Substitute (3.6) and (3.7) into (3.5) to get following relationship:

$$\sum_{n=0}^{\infty} u_n(x,y) = \sum_{k=0}^{m-1} \frac{\partial^k}{\partial y^k}u(x,0^+)\frac{y^k}{k!} + J_y^\alpha(g(x)) - J_y^\alpha(A_n) \qquad (3.11)$$

This equation leads to the recursive relationship by the following iterates:

$$u_0(x,y) = \sum_{k=0}^{m-1} \frac{\partial^k}{\partial y^k}u(x,0^+)\frac{y^k}{k!} + J_y^\alpha(g(x))$$

$$u_1(x,y) = -J_y^\alpha(A_0)$$

$$u_2(x,y) = -J_y^\alpha(A_1) \qquad (3.12)$$

$$u_3(x,y) = -J_y^\alpha(A_2)$$

$$\vdots$$

$$u_{n+1}(x,y) = -J_y^\alpha(A_n), \ n \geq 0$$

(f) Use the terms of the series $u_0, u_1, u_2, u_3,...$ to get the solution by using the approximation $u(x,y) = \lim_{n \to \infty} \Phi_n$, where $\Phi_n = \sum_{k=0}^{n-1} u_k$. It is valuable to note that the zeroth component $u_0(x,y)$ is defined by all terms that arise from the initial condition. The proceeding terms $u_0(x,y), u_1(x,y), u_2(x,y), ....$ and so on, can be completely determined by the results of previous terms as the recursive relationship is presented among them, so, we can compute all the remaining components $u_n(x,y), n \geq 1$ by the relation given in Eq (3.12).

## 4. Applications

### 4.1. Example 1

Let us consider the fractional partial differential equation in the following form

$$D_y^\alpha u(x,y) + u(x,y) D_x^\beta u(x,y) = x \qquad (4.1)$$

Subject to the initial condition

$$u(x,0) = 1 \qquad (4.2)$$

The exact solution of equation (4.1) is given in [42]

$$u(x, y) = x \tanh(y) + \text{sech}(y) \tag{4.3}$$

Following the ADM formulation presented in section 3, we have,

$$u(x,y) = \sum_{k=0}^{m-1} \frac{\partial^k}{\partial y^k} u(x, 0^+) \frac{y^k}{k!} + J_y^\alpha(x) - J_y^\alpha(Nu) \tag{4.4}$$

Substituting the initial condition (4.2) in (4.4) and using (3.7) to calculate Adomian polynomials $A_n$, yields the following recursive relations:

$$u_0(x,y) = \sum_{k=0}^{m-1} \frac{\partial^k}{\partial y^k} u(x, 0^+) \frac{y^k}{k!} + J_y^\alpha(x)$$

$$u_{n+1}(x,y) = -J_y^\alpha(A_n), \; n \geq 0$$

Using above relationship, the first two terms of the decomposition series are given below

$$u_0(x,y) = 1 + \frac{y^\alpha x}{\Gamma[1+\alpha]},$$

$$u_1(x,y) = -J_y^\alpha(A_0) = -\frac{y^{2\alpha} x^{1-\beta} \left(\frac{1}{\Gamma(2\alpha+1)} + \frac{4^\alpha x \Gamma(\alpha+\frac{1}{2}) y^\alpha}{\sqrt{\pi} \Gamma(\alpha+1)\Gamma(3\alpha+1)}\right)}{\Gamma(2-\beta)},$$

The solution in series form can be written as:

$$u(x,y) = 1 + \frac{y^\alpha x}{\Gamma[1+\alpha]} - \frac{y^{2\alpha} x^{1-\beta} \left(\frac{1}{\Gamma(2\alpha+1)} + \frac{4^\alpha x \Gamma\left(\alpha+\frac{1}{2}\right) y^\alpha}{\sqrt{\pi} \Gamma(\alpha+1)\Gamma(3\alpha+1)}\right)}{\Gamma(2-\beta)} + \cdots \tag{4.5}$$

By putting the values of $\alpha = 1$, $\beta = 1$ in equation (4.5), we can check the accuracy of the approximate solution of equation 4.1 obtained by ADM formulation. Table 1 shows the approximate solutions of examples 1 of fractional order partial differential equation for different values of y = 0.01, 0.05 and 0.1 and x = [0, 1] at different values of $\alpha$ = 0.5, 0.75, and 1 and $\beta$ = 0.5, 0.75, and 1. The absolute error is calculated between approximate and the exact solution by [$u_{Exact}$-$u_{ADM}$] at $\alpha = 1$, $\beta = 1$ for fractional order partial differential equation of the form (4.1). Figure 1 also demonstrates the same results with clear visualization that the approximate solution by extended formulation at $\alpha = 1$, $\beta = 1$ is in excellent agreement with the exact solutions.

*Table 1. Approximate solution of Example 1 for different values of α, β and absolute error at α=1 ,β=1*

| | Approximate Solutions by ADM | Exact | Error |
| --- | --- | --- | --- |

| y | x | α = 0.5, β = 0.5 | α = 0.75, β = 0.75 | α = 1, β = 1 | α = 1, β = 1 | $u_{Exact} - u_{ADM}$ |
|---|---|---|---|---|---|---|
| 0.01 | 0.3 | 1.02826 | 1.00972 | 1.00295 | 1.00295 | $1.78455 \times 10^{-17}$ |
|  | 0.6 | 1.05931 | 1.01991 | 1.00595 | 1.00595 | $1.79697 \times 10^{-17}$ |
|  | 0.9 | 1.09087 | 1.03015 | 1.00895 | 1.00895 | $1.81007 \times 10^{-17}$ |
| 0.05 | 0.3 | 1.05085 | 1.02803 | 1.01374 | 1.01374 | $1.35317 \times 10^{-12}$ |
|  | 0.6 | 1.11205 | 1.06088 | 1.02873 | 1.02873 | $1.36598 \times 10^{-12}$ |
|  | 0.9 | 1.17514 | 1.0942 | 1.04371 | 1.04371 | $1.37878 \times 10^{-12}$ |
| 0.1 | 0.3 | 1.05979 | 1.04063 | 1.02492 | 1.02492 | $3.4865 \times 10^{-10}$ |
|  | 0.6 | 1.13731 | 1.09334 | 1.05482 | 1.05482 | $3.55184 \times 10^{-10}$ |
|  | 0.9 | 1.21761 | 1.14748 | 1.08472 | 1.08472 | $3.61719 \times 10^{-10}$ |

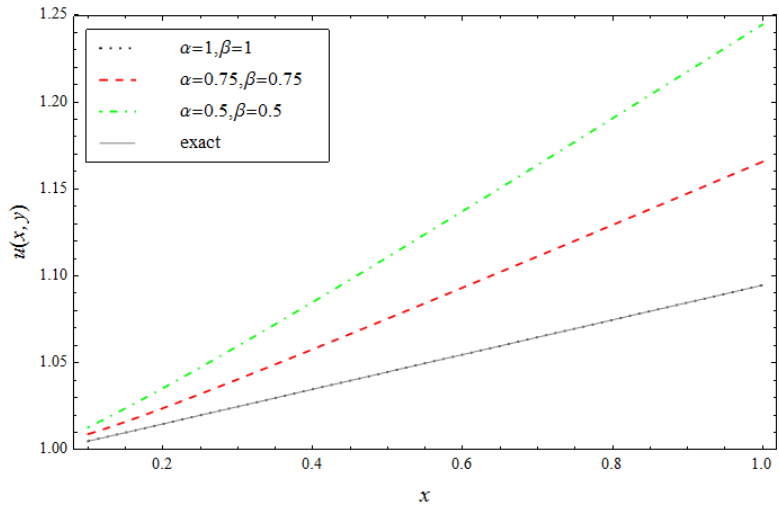

Figure 1. Approximate and exact solution of u(x,y) of Example 1 at different values of α, β and α=1,β=1

### 4.2. Example 2

Let us consider the fractional order partial differential equation in the following form

$$D_y^\alpha u(x,y) + u(x,y) D_x^\beta u(x,y) = 1 \qquad (4.6)$$

Subject to the initial condition

$$u(x,0) = -x \qquad (4.7)$$

The exact solution of equation (4.6) is given in [42]

$$u(x,y) = \frac{2x - 2y + y^2}{2(y-1)} \qquad (4.8)$$

Following the ADM formulation presented in section 3, we have,

$$u(x,y) = \sum_{k=0}^{m-1} \frac{\partial^k}{\partial y^k} u(x,0^+) \frac{y^k}{k!} + J_y^\alpha(1) - J_y^\alpha(Nu) \tag{4.9}$$

Substituting the initial condition (4.7) in (4.9) and using (3.7) to calculate Adomian polynomials $A_n$, yields the following recursive relations:

$$u_0(x,y) = \sum_{k=0}^{m-1} \frac{\partial^k}{\partial y^k} u(x,0^+) \frac{y^k}{k!} + J_y^\alpha(1)$$

$$u_{n+1}(x,y) = -J_y^\alpha(A_n), \; n \geq 0$$

Using above relationship, the first few terms of the decomposition series are given by

$$u_0(x,y) = -x + \frac{y^\alpha}{\Gamma(\alpha+1)},$$

$$u_1(x,y) = -J_y^\alpha(A_0) = \frac{x^{1-\beta} y^\alpha (\Gamma(\alpha+1)y^\alpha - x\Gamma(2\alpha+1))}{\Gamma(\alpha+1)\Gamma(2\alpha+1)\Gamma(2-\beta)},$$

$$u_2(x,y) = -J_y^\alpha(A_1) = x^{1-2\beta} y^{2\alpha} \left( \frac{x^2 \left( \frac{\beta-2}{\Gamma(3-2\beta)} - \frac{1}{\Gamma(2-\beta)^2} \right)}{\Gamma(2\alpha+1)} + \frac{xy^\alpha \left( -\frac{4^\alpha(\beta-2)\Gamma(\alpha+\frac{1}{2})}{\sqrt{\pi}\Gamma(\alpha+1)\Gamma(3-2\beta)} + \frac{1}{\Gamma(2-\beta)^2} + \frac{1}{\Gamma(2-2\beta)} \right)}{\Gamma(3\alpha+1)} - \frac{\Gamma(3\alpha+1)y^{2\alpha}}{\Gamma(\alpha+1)\Gamma(2\alpha+1)\Gamma(4\alpha+1)\Gamma(2-2\beta)} \right),$$

The solution in series form can be written as:

$$u(x,y) = -x + \frac{y^\alpha}{\Gamma(\alpha+1)} + \frac{x^{1-\beta} y^\alpha (\Gamma(\alpha+1)y^\alpha - x\Gamma(2\alpha+1))}{\Gamma(\alpha+1)\Gamma(2\alpha+1)\Gamma(2-\beta)} + + \cdots \tag{4.10}$$

Eq (4.10) gives the highly approximate solution of Eq (4.6). The accuracy of the solution obtained by applying ADM can be evaluated by putting the values of $\alpha = 1$, $\beta = 1$ in equation (4.10). Following Table 2 displays the approximate solutions of examples 2 of fractional order partial differential equations for different values of y = 0.01, 0.05 and 0.1 for x= [0, 1] at different values of $\alpha$ = 0.5, 0.75, and 1 and $\beta$ = 0.5, 0.75, and 1. The absolute error is calculated between approximate and the exact solution by [$u_{Exact}$-$u_{ADM}$] at $\alpha = 1$, $\beta = 1$ for fractional order partial differential equation of the form (4.6). Figure 2 also gives the clear picture of the approximate and exact solutions of example 2 at $\alpha = 1$, $\beta = 1$. It can be easily observed that both the solutions are in line with each other with high accuracy.

*Table 2. Approximate solution of Example 2 for different values of α, β and absolute error at α=1, β=1*

| | | Approximate Solutions by ADM | | | Exact | Error |
|---|---|---|---|---|---|---|
| y | x | $\alpha = 0.5,$ $\beta = 0.5$ | $\alpha = 0.75,$ $\beta = 0.75$ | $\alpha = 1, \beta = 1$ | $\alpha = 1, \beta = 1$ | $u_{Exact}$-$u_{ADM}$ |
| 0.01 | 0.3 | -0.210064 | -0.274905 | -0.29298 | -0.29298 | $2.9798 \times 10^{-9}$ |
|      | 0.6 | -0.555538  | -0.586408  | -0.59601  | -0.59601  | $6.0101 \times 10^{-9}$  |
|      | 0.9 | -0.910206  | -0.898583  | -0.89904  | -0.89904  | $9.04041 \times 10^{-9}$ |
| 0.05 | 0.3 | -0.0782081 | -0.213214  | -0.264472 | -0.264474 | $1.80921 \times 10^{-6}$ |
|      | 0.6 | -0.50313   | -0.556151  | -0.580259 | -0.580263 | $3.78289 \times 10^{-6}$ |
|      | 0.9 | -0.966632  | -0.90218   | -0.896047 | -0.896053 | $5.75658 \times 10^{-6}$ |
| 0.1  | 0.3 | 0.0446003  | -0.147862  | -0.22775  | -0.227778 | $2.77778 \times 10^{-5}$ |
|      | 0.6 | -0.454211  | -0.528458  | -0.56105  | -0.561111 | $6.11111 \times 10^{-5}$ |
|      | 0.9 | -1.03523   | -0.916113  | -0.89435  | -0.894444 | $9.44444 \times 10^{-5}$ |

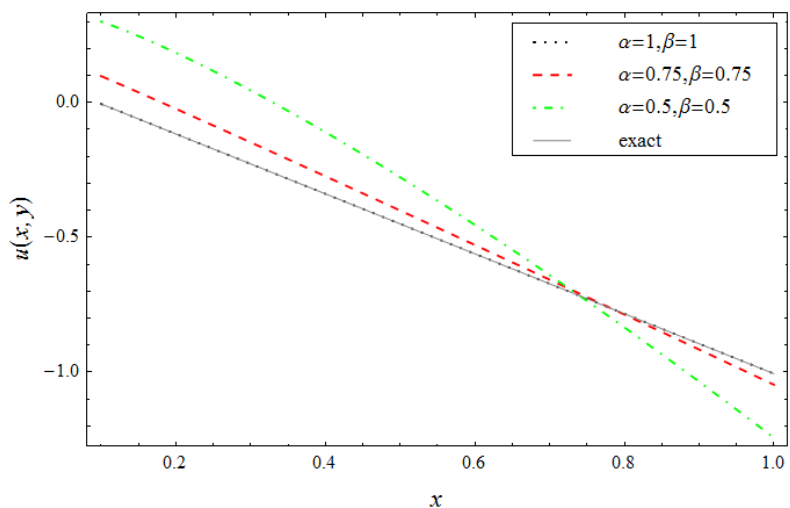

Figure 2. Approximate and exact solution of u(x,y) of Example 2 at different values of α, β and α=1, β=1

### 4.3. Example 3

Let us consider the fractional order partial differential equation in the following form having $g(x) = 0$,

$$D_y^\alpha u(x,y) + u(x,y) D_x^\beta u(x,y) = 0 \tag{4.11}$$

Subject to the initial condition

$$u(x,0) = x + 1 \tag{4.12}$$

The exact solution of equation (4.11) is given in [42]

$$u(x,y) = \frac{1+x}{1+y} \tag{4.13}$$

Following the ADM formulation presented in section 3, we have,

$$u(x,y) = \sum_{k=0}^{m-1} \frac{\partial^k}{\partial y^k} u(x, 0^+) \frac{y^k}{k!} - J_y^\alpha(Nu) \qquad (4.14)$$

Substituting the initial condition (4.12) in (4.14) and using (3.7) to calculate Adomian polynomials $A_n$, yields the following recursive relations:

$$u_0(x,y) = \sum_{k=0}^{m-1} \frac{\partial^k}{\partial y^k} u(x, 0^+) \frac{y^k}{k!}$$

$$u_{n+1}(x,y) = -J_y^\alpha(A_n), \; n \geq 0$$

Using above relationship, the first few terms of the decomposition series are given by

$$u_0(x,y) = x + 1,$$

$$u_1(x,y) = -J_y^\alpha(A_0) = -\frac{(x+1)y^\alpha x^{1-\beta}}{\Gamma(\alpha+1)\Gamma(2-\beta)},$$

$$u_2(x,y) = -J_y^\alpha(A_1) = \frac{(x+1)y^{2\alpha} x^{1-2\beta}\left(\frac{x}{\Gamma(2-\beta)^2} + \frac{2(x+1)-\beta(x+2)}{\Gamma(3-2\beta)}\right)}{\Gamma(2\alpha+1)},$$

By adding these terms, we can get the solution in series form:

$$u(x,y) = x + 1 - \frac{(x+1)y^\alpha x^{1-\beta}}{\Gamma(\alpha+1)\Gamma(2-\beta)} + \frac{(x+1)y^{2\alpha} x^{1-2\beta}\left(\frac{x}{\Gamma(2-\beta)^2} + \frac{2(x+1)-\beta(x+2)}{\Gamma(3-2\beta)}\right)}{\Gamma(2\alpha+1)} + \cdots \qquad (4.15)$$

The approximate solutions can be obtained by putting different values of $\alpha$, $\beta$ in Eq (4.15). Table 3 shows the accuracy of the approximate solution of example 3, calculated by putting $\alpha = 1$, $\beta = 1$ in (4.15). Table 3 also exhibits the approximate solutions of examples 3 of fractional order partial differential equation for different values of y = 0.01, 0.05 and 0.1 and x= [0, 1] at different values of $\alpha$ = 0.5, 0.75, and 1 and $\beta$ = 0.5, 0.75, and 1. The same results of approximate solutions at different values of $\alpha$, $\beta$ and exact solution can be visualized in Figure 3.

Table 3. Approximate solution of Example 3 for different values of $\alpha$, $\beta$ and absolute error at $\alpha$=1, $\beta$=1

| y | x | Approximate Solutions by ADM | | | Exact | Error |
|---|---|---|---|---|---|---|
| | | $\alpha = 0.5,$ $\beta = 0.5$ | $\alpha = 0.75,$ $\beta = 0.75$ | $\alpha = 1, \beta = 1$ | $\alpha = 1, \beta = 1$ | $u_{Exact}-u_{ADM}$ |
| 0.01 | 0.3 | 1.20487 | 1.26054 | 1.28713 | 1.28713 | $1.28713 \times 10^{-8}$ |
| | 0.6 | 1.45717 | 1.54776 | 1.58416 | 1.58416 | $1.58416 \times 10^{-8}$ |
| | 0.9 | 1.71169 | 1.83537 | 1.88119 | 1.88119 | $1.88119 \times 10^{-8}$ |
| 0.05 | 0.3 | 1.10925 | 1.1828 | 1.23809 | 1.2381 | $7.7381 \times 10^{-6}$ |
| | 0.6 | 1.29774 | 1.4429 | 1.5238 | 1.52381 | $9.52381 \times 10^{-6}$ |

| | 0.9 | 1.49262 | 1.70524 | 1.80951 | 1.80952 | $1.13095 \times 10^{-5}$ |
|---|---|---|---|---|---|---|
| 0.1 | 0.3 | 1.00627 | 1.12089 | 1.1817 | 1.18182 | $1.18182 \times 10^{-4}$ |
| | 0.6 | 1.11627 | 1.35561 | 1.4544 | 1.45455 | $1.45455 \times 10^{-4}$ |
| | 0.9 | 1.23329 | 1.59591 | 1.7271 | 1.72727 | $1.72727 \times 10^{-4}$ |

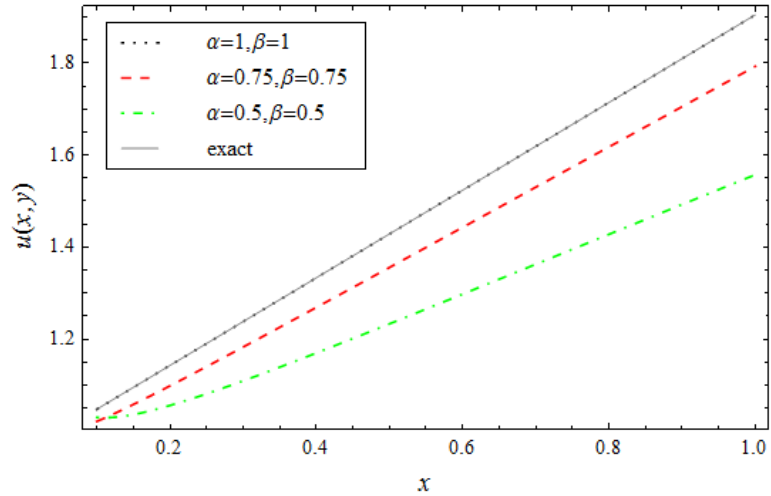

Figure 3. Approximate and exact solution of u(x,y) of Example 3 at different values of α, β and α=1,β=1

### 4.4. Example 4

Let us consider the fractional order partial differential equation in the following form with $g(x) = 0$,

$$D_y^\alpha u(x,y) + u(x,y) D_x^\beta u(x,y) = 0 \qquad (4.16)$$

Subject to the initial condition

$$u(x,0) = x \qquad (4.17)$$

The exact solution of equation (4.16) is given in [42]

$$u(x,y) = \frac{x}{1+y} \qquad (4.18)$$

Following the ADM formulation presented in section 3, we have,

$$u(x,y) = \sum_{k=0}^{m-1} \frac{\partial^k}{\partial y^k} u(x,0^+) \frac{y^k}{k!} - J_y^\alpha(Nu) \qquad (4.19)$$

Substituting the initial condition (4.17) in (4.19) and using (3.7) to calculate Adomian polynomials $A_n$, yields the following recursive relations:

$$u_0(x,y) = \sum_{k=0}^{m-1} \frac{\partial^k}{\partial y^k} u(x,0^+) \frac{y^k}{k!}$$

$$u_{n+1}(x,y) = -J_y^\alpha(A_n), \quad n \geq 0$$

Using above relationship, the first few terms of the decomposition series are given by

$$u_0(x,y) = x,$$

$$u_1(x,y) = -J_y^\alpha(A_0) = -\frac{y^\alpha x^{2-\beta}}{\Gamma(\alpha+1)\Gamma(2-\beta)},$$

$$u_2(x,y) = -J_y^\alpha(A_1) = \frac{\left(\frac{2-\beta}{\Gamma(3-2\beta)} + \frac{1}{\Gamma(2-\beta)^2}\right) y^{2\alpha} x^{3-2\beta}}{\Gamma(2\alpha+1)},$$

The solution in series form can be written as:

$$u(x,y) = x - \frac{y^\alpha x^{2-\beta}}{\Gamma(\alpha+1)\Gamma(2-\beta)} + \frac{\left(\frac{2-\beta}{\Gamma(3-2\beta)} + \frac{1}{\Gamma(2-\beta)^2}\right) y^{2\alpha} x^{3-2\beta}}{\Gamma(2\alpha+1)} + \cdots \quad (4.20)$$

Adomian decomposition method is a very reliable and efficient tool that gives the highly accurate solution. Like previous examples (1-3), example 4 also provides the extremely similar approximate solution as exact solution. It is also evaluated by calculating absolute error between approximate and the exact solution by [$u_{Exact}$-$u_{ADM}$] at $\alpha = 1$, $\beta = 1$ for fractional order partial differential equation of the form (4.16). The results are listed in Table 4. Figure 4 also reveals that the pattern of the approximate solution by extended ADM formulation at $\alpha = 1$, $\beta = 1$ promises the exact solutions. It is also perceived that as we move along the domain, we get consistent accuracy.

Table 4. Approximate solution of Example 4 for different values of α, β and absolute error at α=1, β=1

| y | x | Approximate Solutions by ADM | | | Exact | Error |
|---|---|---|---|---|---|---|
| | | $\alpha = 0.5, \beta = 0.5$ | $\alpha = 0.75, \beta = 0.75$ | $\alpha = 1, \beta = 1$ | $\alpha = 1, \beta = 1$ | $u_{Exact}$-$u_{ADM}$ |
| 0.01 | 0.3 | 0.276009 | 0.290771 | 0.29703 | 0.29703 | $2.97029 \times 10^{-13}$ |
| | 0.6 | 0.544279 | 0.580275 | 0.594059 | 0.594059 | $5.94058 \times 10^{-13}$ |
| | 0.9 | 0.80891 | 0.869243 | 0.891089 | 0.891089 | $8.91087 \times 10^{-13}$ |
| 0.05 | 0.3 | 0.252999 | 0.271796 | 0.285714 | 0.285714 | $4.46429 \times 10^{-9}$ |
| | 0.6 | 0.491149 | 0.540065 | 0.571429 | 0.571429 | $8.92857 \times 10^{-9}$ |
| | 0.9 | 0.720922 | 0.806873 | 0.857143 | 0.857143 | $1.33929 \times 10^{-8}$ |
| 0.1 | 0.3 | 0.23591 | 0.256139 | 0.272727 | 0.272727 | $2.72727 \times 10^{-7}$ |
| | 0.6 | 0.442692 | 0.507181 | 0.545454 | 0.545454 | $5.45455 \times 10^{-7}$ |
| | 0.9 | 0.624414 | 0.756131 | 0.818181 | 0.818181 | $8.18182 \times 10^{-7}$ |

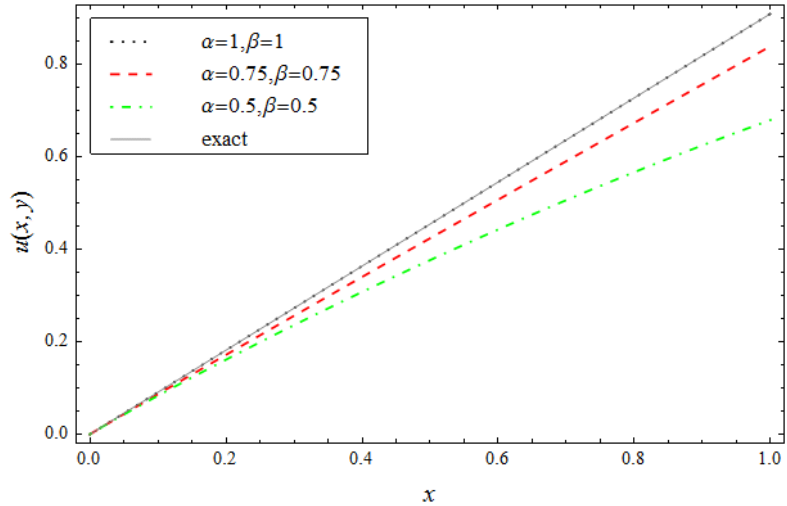

*Figure 4. Approximate and exact solution of u(x,y) of Example 4 at different values of α, β and α=1,β=1*

Remarks:

The solution in series form is a convergent series. The convergence of the series can be further improved by computing more and more terms in the decomposition method i.e. by increasing the number of iterations one can reach desired accuracy.

## 5. Conclusion

The objective of this study is to implement the Adomian decomposition method to obtain explicit and optimal numerical solution of fractional order partial differential equations with initial conditions. The obtained solution demonstrates the reliability of the proposed method and shows the applicability of this method to solve large amount of non-linear fractional equations. The proposed method is free from assumptions and rounding-off errors. It offers more realistic series solutions whose continuity depends on the fractional derivative. Convergence of approximate solutions is observed as the number of decomposed terms is increased. Consistent accuracy throughout the domain of the problem is also witnessed.